\documentclass[11pt]{amsart}

\usepackage{epigamath}

\usepackage[english]{babel}

\usepackage{amscd, amssymb, amsmath, mathrsfs, amsthm}
\usepackage[utf8]{inputenc}
\usepackage{tikz-cd}
\usetikzlibrary{decorations.pathreplacing}
\usepackage{booktabs}
\usepackage{color}

\usepackage{microtype}
\usepackage{enumerate}
\usepackage{enumitem}

\usepackage{xypic}
\usepackage{graphicx}
\usepackage{fancyhdr} % used by RCS
\usepackage{enumerate}
\newtheorem{theorem}[equation]{Theorem}      % (If you want theorem numbered
\newtheorem{lemma}[equation]{Lemma}          %
\newtheorem{corollary}[equation]{Corollary}  %       goes for lemmas, etc.)
\newtheorem{proposition}[equation]{Proposition}

\theoremstyle{definition}
\newtheorem{conjecture}[equation]{Conjecture}

\newtheorem{question}[equation]{Problem}

\theoremstyle{definition}
\newtheorem{defn}[equation]{Definition}
\theoremstyle{remark}
\newtheorem{rem}[equation]{Remark}

\theoremstyle{definition}
\newtheorem{remark}[equation]{Remark}

\numberwithin{equation}{section}
\newcommand{\be}{\begin{equation}}
\newcommand{\ee}{\end{equation}}
\newcommand{\bes}{\begin{equation*}}
\newcommand{\ees}{\end{equation*}}
\newcommand{\bea}{\begin{eqnarray}}
\newcommand{\eea}{\end{eqnarray}}
\newcommand{\beas}{\begin{eqnarray}}
\newcommand{\eeas}{\end{eqnarray}}
\newcommand{\ben}{\begin{note}}
\newcommand{\een}{\end{note}}
\newcommand{\bexl}{\vskip0.1em\noindent\hrulefill\vskip1em\begin{ExerciseList}}
\newcommand{\eexl}{\end{ExerciseList}\hrulefill}

\newcommand{\bthm}{\begin{theorem}}
\newcommand{\ethm}{\end{theorem}}
\newcommand{\bpro}{\begin{prop}}
\newcommand{\epro}{\end{prop}}
\newcommand{\bcor}{\begin{corollary}}
\newcommand{\ecor}{\end{corollary}}
\newcommand{\bcon}{\begin{conjecture}}
\newcommand{\econ}{\end{conjecture}}
\newcommand{\bp}{\begin{proof}}
\newcommand{\ep}{\end{proof}}
\newcommand{\blem}{\begin{lemma}}
\newcommand{\elem}{\end{lemma}}
\newcommand{\bn}{\begin{note}}
\newcommand{\en}{\end{note}}
\newcommand{\benum}{\begin{enumerate}}
\newcommand{\eenum}{\end{enumerate}}
\newcommand{\bed}{\begin{defn}}
\newcommand{\eed}{\end{defn}}
\newcommand{\brem}{\begin{remark}}
\newcommand{\erem}{\end{remark}}

\newcommand{\btik}{\begin{tikzpicture}\begin{axis}[scale=0.5,axis y line=center, axis x line=middle]}
\newcommand{\etik}{\end{axis}\end{tikzpicture}}

\let\into=\hookrightarrow
\let\mapsto=\longmapsto

\newcommand{\upperRomannumeral}[1]{\uppercase\expandafter{\romannumeral#1}}

\let\into=\hookrightarrow
\let\isom=\simeq
\let\rk=\rank
\let\tensor=\otimes

\newcommand{\End}{\mathscr{E}nd}

\newcommand{\Hom}{{\rm Hom}}

\newcommand{\mydot}{{\scriptstyle{\bullet}}}

\renewcommand{\int}{\operatorname{int}}
\renewcommand{\O}{{\mathcal O}}

\renewcommand{\bpro}{\begin{proposition}}
\renewcommand{\epro}{\end{proposition}}

\newcommand{\edim}{{\rm exp. \dim}}
\newcommand{\fs}[1]{\mathit{F}_*({#1})}
\newcommand{\shom}{\mathscr{H}om}
\newcommand{\qquot}{{{\rm Quot}^{r,0}(\fs{Q})}}
\newcommand{\rk}{{\rm rk}}
\newcommand{\trip}[1]{({#1},\nabla,{#1}_{\mydot})}

\newcommand{\Gr}[1]{\mathscr{G}_{#1}}
\newcommand{\qqquot}[1]{{{\rm Quot}^{r,0}(\fs{#1})}}
\newcommand{\sC}{\mathscr{C}}
\newcommand{\codim}{{\rm codim}}

\EpigaVolumeYear{4}{2020}
\EpigaArticleNr{17}
\ReceivedOn{August 28, 2019}
\InFinalFormOn{June 30, 2020}
\AcceptedOn{September 13, 2020}

\title{Opers of higher types, Quot-schemes and\\Frobenius instability loci}
\titlemark{Opers of higher types, Quot-schemes and Frobenius instability loci}

\author{Kirti Joshi}
\address{Mathematics Department, University of Arizona, 617 N Santa Rita, Tucson
85721-0089, USA} \email{kirti@math.arizona.edu}

\author{Christian Pauly}
\address{Laboratoire de Math\'ematiques J.A. Dieudonn\'e, UMR 7351 CNRS, Universit\'e de Nice Sophia-Antipolis, 06108 Nice Cedex 02, France} 
\email{pauly@unice.fr}

\authormark{K. Joshi and C. Pauly}

\AbstractInEnglish{\small{In this paper we continue our study of the Frobenius instability locus
in the coarse moduli space of semi-stable vector bundles of rank $r$ and degree $0$
over a smooth projective curve defined over an algebraically closed field of
characteristic $p>0$. In a previous paper we identified the ``maximal"
Frobenius instability strata with opers (more precisely as opers of type $1$ in the terminology of the present paper) and related them
to certain Quot-schemes of Frobenius direct images of line bundles. The main
aim of this paper is to describe for any integer $q \geq 1$ a conjectural
generalization of this correspondence between opers of type $q$ and
Quot-schemes of Frobenius direct images of vector bundles of rank $q$.
We also give a conjectural formula for the dimension of the Frobenius 
instability locus.}}



\TitleInFrench{Opers de type supérieur, schémas Quot et lieux d'instabilité de Frobenius}

\AbstractInFrench{\small{Nous continuons l'étude du lieu d'instabilité de Frobenius dans l'espace de modules grossier des fibrés semi-stables de rang $r$ et degré $0$ sur une courbe projective lisse définie sur un corps algébriquement clos de caractéristique $p>0$. Dans un article précédent, nous avons identifié les strates «~maximales~» d'instabilité de Frobenius avec certains opers (de type $1$, avec la terminologie du présent article) et avons montré leur lien avec certains schémas Quot d'images directes sous Frobenius de fibrés en droites. Le principal objectif de cet article est de décrire, pour tout $q\geq 1$, une généralisation conjecturale de cette correspondance entre opers de type $q$ et schémas Quot d'images directes sous Frobenius de fibrés de rang $q$. Nous donnons également une formule conjecturale pour la dimension du lieu d'instabilité de Frobenius.}}

\acknowledgement{KJ would like to thank the University of Nice Sophia-Antipolis for  
financial support of research visits in June 2013 and June 2016.}

\begin{document}

%%%%%%%%%%%%%%%%%%%%%%%%%%%%%%%
% Title page
%%%%%%%%%%%%%%%%%%%%%%%%%%%%%%%

\removeabove{0,8cm}
\removebetween{0,8cm}
\removebelow{0,8cm}

\maketitle

\begin{prelims}

\DisplayAbstractInEnglish

\bigskip

\DisplayKeyWords

\medskip

\DisplayMSCclass

\bigskip

\languagesection{Fran\c{c}ais}

\bigskip

\DisplayTitleInFrench

\medskip

\DisplayAbstractInFrench

\end{prelims}

%%%%%%%%%%%%%%%%%%%%%
% Table of Contents
%%%%%%%%%%%%%%%%%%%%%

\newpage

\setcounter{tocdepth}{1} 

\tableofcontents

%%%%%%%%%%%%%%%%%%%%%
% Content begins here
%%%%%%%%%%%%%%%%%%%%%

\section{Introduction}
Let $p$ be a prime number. Let $k$ be an algebraically closed field of characteristic $p$. Let $X/k$ be a connected, smooth projective curve over $k$. We will write $F:X\to X$ for the absolute Frobenius morphism of $X$. A foundational classical problem in the theory of vector bundles on smooth, projective curves is the following: 
\begin{question} \label{q:1}
Describe the Frobenius instability locus, i.e. the locus of all stable (and also semi-stable) vector bundles $V$ over $X$ such 
that $F^*(V)$ is not  semi-stable.
\end{question}
This goal has been partially achieved by us  in \cite{joshi15} where, for $p\gg 0$ we provided an explicit construction of all stable bundles $V$ such that $F^*(V)$ is not semi-stable. By a well-known theorem of \cite{shatz76}, one may equip the moduli space of semi-stable vector bundles of fixed degree and rank by a stratification defined using Harder-Narasimhan polygons of $F^*(V)$.
The main problem addressed in the present paper is the finer problem:
\begin{question} 
Describe all Frobenius instability strata.
\end{question}
One of the important results of \cite{joshi15} asserts that there is an explicit polygon, called the oper polygon (which is always attained in every genus and rank) with the property that the Harder-Narasimhan polygon of every vector bundle $F^*(V)$, with $V$ stable, 
lies on or strictly below the oper polygon. 

So to clarify the second problem, one can for example ask: do all possible Harder-Narasimhan polygons occur?
This seems substantially more difficult and complete results are presently available only for $p=2$ (see \cite{joshi-xia00}). 
We have obtained a number of results in small rank, genus, and characteristic, which we will report in a companion paper under preparation. Recently \cite{li16} has also studied the genus-$2$, rank-$3$, characteristic-$3$ case in detail (our small genus, rank, characteristic results include his results, but by methods quite different from his).

The description which we provided in our answer to Problem \ref{q:1} is that all Frobenius-destabilized stable bundles arise 
from suitable quotients of bundles of the form $F_*(Q)$ with $Q$ a stable vector bundle. In all the small genus, small rank, small characteristic situations, the vector bundle $Q$ is of rank $1$.  If $\mathrm{rk}(Q) =1$, the Quot-schemes $\qquot$ 
parameterizing subsheaves of given rank $r$ and degree $0$ of $F_*(Q)$ are substantially better behaved than for general rank :
in fact, one of the main results of \cite{joshi15} says that if $\mathrm{rk}(Q) =1$ then the dimension of
$\qquot$ is $0$ if its expected dimension is $0$, and its $k$-rational points correspond to opers (of type $1$). As shown
in section 6, this fails if $\rk(Q) > 1$. The existence of these higher-dimensional components of $\qquot$ stems from
certain line subbundles $L \hookrightarrow Q$ of sufficiently high degree, which induce natural embeddings
$\mathrm{Quot}^{r,0}(F_*(L)) \hookrightarrow \qquot$.

In the present paper we lay the ground work for addressing the second problem and its finer special cases. The key tool 
we introduce here is the notion of opers of type $q$ and rank $\ell q$ (see \cite[Definition~3.1.1]{joshi15}, recalled in Section 4). This notion generalizes the notion of an oper
$(V, \nabla, V_\mydot)$ studied by \cite{beilinson00b} and consisting of a vector bundle $V$ equipped
with an integrable connection $\nabla$ and a full flag $V_\mydot$ satisfying some transversality conditions. An oper of 
type $q$ should be thought of as the bundle analog of the parabolic induction (and its adjoint the Jacquet Functor) in 
the theory of automorphic forms --- the parabolic in the present case is of type given by the $\ell$-tuple $(q,\ldots,q)$, 
which corresponds to a partial flag whose associated quotients are all of rank $q$.
Given an oper of type $q$ and rank $\ell q$ we naturally obtain a rank-$q$ bundle $Q$ as first quotient of the oper filtration.

In the case $\rk(Q) >1$ we conjecture (Conjecture 7.6) a similar statement assuming that one restricts attention
to ``non-degenerate" subsheaves of $F_*(Q)$, i.e. excluding in particular the above-mentioned sub-Quot-schemes
$\mathrm{Quot}^{r,0}(F_*(L))$. More precisely, we conjecture that, if the expected dimension of $\mathrm{Quot}^{\ell q,0}(F_*(Q))$
is $0$, then dormant opers of type $q = \rk(Q)$ form a non-empty open subset of dimension $0$ of the 
Quot-scheme $\mathrm{Quot}^{\ell q,0}(F_*(Q))$, which has, as mentioned above, components of dimension $>0$.
Here dormant means that the $p$-curvature of the connection $\nabla$ is zero.
We check that this conjecture holds when $Q$ is a semi-stable direct sum of $q$
line bundles (Theorem 3.5).

Finally, somewhat independently of the previous considerations, we give a conjecture for the dimension of the
Frobenius instability locus $\mathcal{J}(r)$ in the coarse moduli space of semi-stable rank $r$ and degree $0$
vector bundles. This conjecture says that 
$$\dim \mathcal{J}(r) \geq (r^2 -r +1)(g-1) - (r-1).$$
The conjecture holds for $r=2$. We also conjecture that a general vector bundle in $\mathcal{J}(r)$
has ``minimal" Harder-Narasimhan filtration (Conjecture 8.3).

\bigskip

\noindent{\bf Acknowledgements.} We would like to thank the referee for useful comments and suggestions.

\section{Preliminaries on Quot-schemes}

In what follows, the following notations and assumptions will be in force. 
Let $X$ be a smooth, projective curve of genus $g\geq 2$  over an algebraically closed field $k$  of
characteristic $p>0$. Let $F: X\to X$ be the absolute Frobenius morphism of $X$. For a vector bundle $V$,
we shall write $V^*$ for its dual $\shom( V, \mathcal{O}_X)$.

\bigskip
Suppose $\theta$ is a theta-characteristic for $X$. Let $\Gr{2}$ be the unique non-split extension of  $\theta^{-1}$ by $\theta$: this bundle depends on the choice of the theta-characteristic $\theta$ but our notation suppresses this dependence. For $r<p$ let  $\Gr{r}={\rm Sym}^{r-1}(\Gr{2})$ which will be called the Gunning bundle of rank $r$. This is an indecomposable bundle of degree zero and trivial determinant.

\bigskip

We recall two formulas from \cite{jrxy06}. Let $Q$ be a vector bundle of rank $q$ and slope $\mu(Q)$ on $X$. Then
\be\label{eq:degree-formula}  
\deg(\fs Q)=q\mu(Q)+q(p-1)(g-1),
\ee
and equivalently
\be\label{eq:slope-formula} 
\mu(\fs Q)=\frac{\mu(Q)}{p}+\frac{(p-1)(g-1)}{p}.
\ee

\bigskip

The first is a consequence of the Riemann-Roch formula $\chi(V)=\deg(V)+\rk(V)(1-g)$ for a vector bundle $V$ of rank $\rk(V)$ and degree $\deg(V)$ on $X$ 
and the fact that $\chi(V) = \chi(\fs V)$ for the finite map $F$. The second is equivalent to the first.

\bigskip

Let $Q$ be a stable bundle on $X$ of rank $q$ and slope $\mu(Q)=\mu<0$. Let $r\geq 1$ be an integer. 
For a coherent sheaf $V$ and a surjection $V \to G$  to a coherent sheaf, we will say that $G$ has codegree $d$ if the kernel $\ker(V \to G)$ has degree $d$. Similarly we will say that $G$ has corank $r$ if $\ker(V \to G)$ has rank $r$. Let $\qquot$ be the Quot-scheme of quotients of $\fs{Q}$ of codegree $0$ and corank $r$. If $V \to G$ is a quotient with kernel $E$ we will 
habitually write $[E]$ for the corresponding point of the relevant Quot-scheme corresponding to this quotient. 

Let $[E]\in\qquot$ be a point of the Quot-scheme. Then we define the integer 
\be 
e(E)=\edim_{[E]}(\qquot)=\chi(\shom(E,\fs{Q}/E)).
\ee
The integer $e(E)$ is called the \textit{expected dimension} of 
$\qquot$ at the point $[E]$ (see \cite[Chapter 2]{huybrechts-lehn-book}).
\bpro\label{p:expected-dim}
Let $r\geq 1$ be an integer. Let $Q$ be a vector bundle on $X$ of rank $q\leq r-1$ and 
slope $\mu(Q)=\mu$. Let $[E]\in \qquot$ be a point of the Quot-scheme. Then the following assertions hold.
\benum[label={\bf(\arabic*)}]
\item The expected dimension of $\qquot$ at $[E]$ is 
$$e(E)= r \deg(Q)+(r^2-q r)(g-1).$$
\item We have $e(E)=0$ if and only if 
$$\deg(Q)=(r-q)(1-g).$$
\item In particular, if $r=\ell q$ for some integer $\ell \geq 2$, then $e(E)=0$ if and only if 
$$\deg(Q)=q(\ell-1)(1-g).$$
\item Moreover if $\mu=-\frac{1}{q}$ (equivalently $\deg(Q)=-1$) then
$$e(E)=(r^2-qr)(g-1)-r.$$
\item If $\deg(Q)=(r-q)(1-g)+d$ with $d\geq 0$, then
$$e(E)=rd.$$
\eenum
\epro 
\bp 
It is enough to prove the first assertion as the rest are immediate consequences of the first.
Let $[E]\in\qquot$ be a point corresponding to a quotient $\fs Q\to G$. Then by definition the expected dimension is given by $e(E)=\chi(\shom(E,G))$. By the Riemann-Roch formula 
we have
$$\chi(\shom(E,G))=\deg(E^*\tensor G)+\rk(E) \rk(G)(1-g)$$
which gives
$$
e(E)=r\deg(G)+r(pq-r)(1-g).
$$
As $E$ has degree zero, so $\deg(G)=\deg(\fs Q)=q\mu+q(p-1)(g-1)$. Substituting this in the above equation and simplifying the result gives the asserted formula. 
\ep

\bpro \label{dimensionirrcom}
Let $Q$ be a semistable bundle of rank $q$ and $\deg(Q)=(r-q)(1-g)+d$ with $q<r<pq$ and $d\geq 0$.
Then
\benum[label={\bf(\arabic*)}]
	\item The Quot-scheme $\qqquot{Q}$ is non-empty.
	\item Any irreducible component of $\qqquot{Q}$ has dimension $\geq rd$.
	\item If $p>\max(2r(r-1)(g-1),\frac{2rd}{q})$ then for any $[E]\in\qqquot{Q}$, the quotient $\fs{Q}/E$ is torsion-free, \textit{i.e.,} $E\into \fs{Q}$ is a subbundle.
\eenum
\epro
\bp 
The assertion (1) is \cite[Proposition 2.3.2]{joshi15}. The assertion (2) is 
\cite[Proposition 2.3.4]{joshi15}. So it remains to prove (3).
Let $E'$ denote the saturation of $E$ in $\fs{Q}$. Then $\mu(E') \geq 0$. 
Assume on the contrary that $E' \not= E$, i.e., $\mu(E')>0$. Then we invoke \cite[Proposition 4.2.1]{joshi15} with  
$n=r$ and $\delta=\frac{1}{2r}$, and we note that $\mu(Q)=\frac{(r-q)(1-g)}{q}+\frac{d}{q}$ and that we have the following
inequality
$$\frac{\mu(Q)}{p}+\delta=\frac{(r-q)(1-g)}{pq}+\frac{d}{pq}+\frac{1}{2r}<\frac{1}{2r} + \frac{1}{2r} = \frac{1}{r}.$$  
Then \cite[Proposition 4.2.1]{joshi15} implies that
$$ 0 < \mu(E') < \frac{\mu(Q)}{p}+\delta < \frac{1}{r}, $$
which leads to a contradiction, since $\rk(E') = \rk(E) = r$.
\ep

\bigskip

\section{A finiteness theorem}

\bigskip

We start with a lemma.

\newcommand{\qqq}{{\rm Quot}^{r',0}(\fs{Q})}
\blem\label{p:empty-quot}
	Let $Q$ be a line bundle of degree $-(r-1)(g-1)$ with $r\geq 2$. If $r'<r$ then 
	$$\qqq=\emptyset.$$
\elem

\bigskip

\bp 
Suppose the assertion is not true. Then we have $\qqq\neq\emptyset$. Consider a subsheaf $E\in\qqq$. We will use the notation of \cite[Lemma 3.4.2]{joshi15}. 
Let $W=F^*(E)$ and equip it with the filtration induced by the
natural filtration $V_{\bullet}$ on the bundle $V = F^*(\fs{Q})$. We shall denote the latter by
$0 = W_{m+1} \subset W_m\subseteq\cdots\subseteq W_1\subseteq W_0=W$. Let $r_i=\rk(W_i/W_{i+1})$. Then $1=r_0\geq r_1\geq \cdots \geq r_m\geq 1$. So $r_i=1$ for all $i\geq 0$. Now we have the following inequalities
\beas 
 0 = \deg(W) = \sum_{i=0}^m\deg(W_i/W_{i+1})& \leq & \sum_{i=0}^m\deg(V_i/V_{i+1}), \\
0&\leq&\sum_{i=0}^m(\deg(Q)+i(2g-2)),\\
0&\leq& (m+1)\deg(Q)+\frac{m(m+1)}{2}(2g-2).
\eeas
As $\rk(E)=r '$ so $m=r'-1$ and hence the last inequality can be written as 
$$0\leq r'(r-1)(1-g)+r'(r'-1)(g-1)=r'(g-1)(r'-1-r+1)=r'(r'-r)(g-1)<0.$$
Thus we have arrived at a contradiction.
\ep 

\newcommand{\quot}[2]{{\rm Quot^{#1}(\fs{#2})}}
\newcommand{\sQ}{{\mathscr Q}}
\newcommand{\sF}{{\mathscr F}}
\newcommand{\sG}{{\mathscr G}}
\newcommand{\sV}{{\mathscr V}}
\newcommand{\sW}{{\mathscr W}}

We will now prove a finiteness theorem for the Quot-scheme $\quot{\ell q,0}{Q}$ when $Q$ is a 
decomposable bundle.
Assume that $q \geq 2, \ell \geq 2$ and that $p > r(r-1)(r-2)(g-1)$
for $r = \ell q$.
Let $L_1,\ldots,L_q$ be $q$  distinct line bundles of degree $-(\ell-1)(g-1)$. We denote by $Q$ the decomposable bundle 
$$Q = \oplus_{i=1}^q L_i \ \ \text{and} \ \  Q_j = \oplus_{i \not= j} L_i \subset Q.$$ 
Then for any $j = 1, \ldots, q$ we have an
inclusion
$$\quot{\ell q,0}{Q_j} \hookrightarrow \quot{\ell q,0}{Q}.$$
By Proposition \ref{dimensionirrcom}(2) we have 
$$\dim \quot{\ell q,0}{Q_j} \geq q \ell (\ell -1)(g-1) > 0.$$ 
We denote by
$\Omega(Q)$ the residual component of these $q$ Quot-schemes in $\quot{\ell q,0}{Q}$, i.e., we have a union
$$ \quot{\ell q,0}{Q} = \bigcup_{j=1}^q \quot{\ell q,0}{Q_j}  \cup \Omega(Q). $$

\bthm\label{p:quot-product-lemma} 
We consider the morphism
$$ \Phi: \prod_{i=1}^q \quot{\ell,0}{L_i} \longrightarrow \quot{\ell q,0}{Q} $$
defined by
$$(S_i \subset \fs{L_i})_{i=1}^q \mapsto (\oplus_{i=1}^q S_i \subset \fs{Q}). $$
Then $\Phi$ induces a bijection at the level of $k$-rational
points
$$\prod_{i=1}^q\quot{\ell,0}{L_i}(k) \isom \Omega(Q)(k).$$
In particular $\Omega(Q)(k) \not= \emptyset$ and $$\dim \Omega(Q) =0.$$
\ethm

\bigskip

\bp 
Let $[E] \in \quot{\ell q,0}{Q}$. We will denote by $E_i = \pi_i(E)$ the image of $E$ 
under the projection 
$\pi_i : F_*(Q) \rightarrow F_*(L_i)$.  First we note the equivalence $[E] \in \quot{\ell q,0}{Q_i}$ if and only if
$E_i = \pi_i(E) = 0$. Therefore, for any $[E] \in \left(\quot{\ell q,0}{Q} \setminus \bigcup_{i=1}^q \quot{\ell q,0}{Q_i}\right)$ we have $E_i \not= 0$ for $i= 1, \ldots, q$. Moreover, by \cite[Proposition~4.5.1 and Theorem~4.1.1]{joshi15} the vector bundles $E$ and $E_i$ 
are semi-stable and their degrees are zero. We now consider the kernel 
$$K_i=\ker(E\to \oplus_{j\neq i}E_j).$$
Then clearly $K_i\subset \fs{L_i}$. Note that $K_i= E\cap \fs{L_i}$ and so 
$K_i\cap K_j=0$ for $i\neq j$. Moreover, by our assumption on $E$, we have 
$K_i \not= 0$ for all $i= 1, \ldots, q$. As the bundles $E$ and $E_i$ are 
semi-stable of degree zero, the bundles $K_i$ are also semi-stable of degree zero. 
We now apply Lemma \ref{p:empty-quot} to the line bundle $L_i$ of degree $-(\ell -1)(g-1)$ and use the fact
that $\quot{\rk(K_i),0}{L_i} \not= \emptyset$ --- since this Quot-scheme contains $[K_i]$ --- and
we obtain that $\rk(K_i)\geq \ell$. Thus we have a homomorphism
$$\oplus_{i=1}^q K_i\to E$$ 
which is injective and so both bundles have the same rank and degree.  So this map is an isomorphism.
Moreover, we observe that $K_i = E_i$ for any $i = 1, \ldots, q$. Thus we have shown that $\Phi$ induces a 
bijection at the level of $k$-rational points between $\prod_{i=1}^q\quot{\ell,0}{L_i}(k)$ and 
$\Omega(Q)(k)$.

Since by \cite[Theorem 6.2.1]{joshi15} each $\quot{\ell,0}{L_i}(k)$ is a finite set, so $\Omega(Q)(k)$ is also a finite
set. Since $k$ is algebraically closed, we deduce that $$\dim \Omega(Q)=0.$$
\ep

\begin{rem} If the line bundles $L_i$ are not distinct, one still can show finiteness of 
$\Omega(Q)$ defined as the residual component of the Quot-schemes $\quot{\ell q,0}{\tilde{Q}}$ for
any subbundle $\tilde{Q} \subset Q$ of degree $0$ and rank $q-1$.
\end{rem}

\bigskip

\begin{rem}
We will define in Section 7 an open subscheme $\Omega(Q) \subset \quot{\ell q,0}{Q}$  for
any semi-stable vector bundle $Q$, generalizing the above $\Omega(Q)$.
\end{rem}

\section{Opers of type $q$}

We now recall the definition of our main object  (see \cite[Definition 3.1.1]{joshi15}).

\newcommand{\gr}[2]{{\rm gr}^{#1}(#2)}

\begin{defn}
An oper of type $q$ on a smooth projective curve $X$ is a triple $\trip V$, where
\benum[label={\bf(\arabic*)}]
\item $V$ is a vector bundle of rank $r=\ell q$ on $X$,
\item $\nabla:V\to V\tensor\Omega^1_X$ is an integrable connection on $V$,
\item $V_\mydot$ is a filtration of length $\ell$ on $V$
$$0=V_\ell\subset V_{\ell-1}\subset \cdots \subset V_1\subset V_0=V$$
such that
\benum[label={\bf(\alph*)}]
\item $V_0/V_1=Q$ is a vector bundle of rank $q$, and
\item $\nabla$ induces $\mathcal{O}_X$-linear isomorphisms for $1 \leq i \leq \ell-1$
$$\gr{i}{V} \to \gr{i-1}{V} \tensor \Omega^1_X,$$
where we define $\gr{i}{V} := V_i/ V_{i+1}$ for all $0 \leq i \leq \ell-1$.
\eenum

\eenum
\end{defn}

\bigskip

Note that our terminology suppresses the dependence on $Q$, but occasionally we may need to emphasize the dependence on $Q$ and in such contexts will refer to such a triple as an oper of type $Q$.

From the definition it is immediate that one has, for all $0\leq i\leq \ell-1$, isomorphisms of bundles
\be 
\gr{i}{V} \isom Q\tensor(\Omega^1_X)^{\tensor i},
\ee 
and 
\be 
\det(V_i)=\det(Q)^{\ell-i}\tensor (\Omega^1_X)^{q((\ell-1)+(\ell-2)+\cdots +i)}=\det(Q)^{\ell-i}\tensor(\Omega_X^1)^{\frac{q(\ell-i)(\ell+i-1)}{2}}
\ee
Note that we have the following equivalence
$$ \deg (Q) = -q(\ell -1)(g-1) \ \iff \ \deg(V) = 0.$$
In the next proposition we will show that opers of type $q$ exist over an algebraically closed field of any characteristic. 
We first recall that $\Gr \ell$ denotes the Gunning bundle of rank $\ell$ associated to the theta-characteristic $\theta$. 
We equip $\Gr\ell$ with any connection $\nabla_{\Gr{\ell}}$ and with its natural filtration $\Gr\mydot$. Then it is well-known 
that the triple $( \mathscr{G}_\ell  , \nabla_{\Gr{\ell}}, \mathscr{G}_{\mydot} )$
is an oper (of type $1$) under the assumption $p > (\ell-1)(g-1)$, if $\mathrm{char}(k) = p > 0$.

\bpro
Let $k$ be any algebraically closed field. If $\mathrm{char}(k) = p > 0$, we assume that $p > q(\ell-1)(g-1)$. Let $S$ be a stable vector bundle on $X$ with $\deg (S) = 0$ and $\rk (S) = q$.  Then the triple
$$ (V = \mathscr{G}_\ell \otimes S, \nabla , V_\mydot =  \mathscr{G}_\mydot \otimes S ) $$
is an oper of type $q$ for any connection $\nabla$ on $\mathscr{G}_\ell \otimes S$.
\epro

\bp 
The proof is elementary. Let us note that the space of connections on $\mathscr{G}_\ell \otimes S$ is non-empty, since
by stability the degree-$0$ vector bundle $S$ admits a connection $\nabla_S$, so $\nabla_{\mathscr{G}_\ell} \otimes \nabla_S$ is a 
connection on $\mathscr{G}_\ell \otimes S$. Note that, again by stability of $S$, the filtration
$V_\mydot$ coincides with the Harder-Narasimhan filtration of $V$. The assumption on $p$ implies that 
no proper subbundle of the filtration  $\mathscr{G}_\mydot \otimes S$
admits a connection. These observations easily show that any connection $\nabla$ on $V$ is an oper connection. Finally we
note that $V_0/ V_1 = S \otimes \theta^{-(\ell - 1)}$.
\ep

\bigskip

We say that an oper $\trip V$ of type $q$ is \textit{nilpotent} (resp. \textit{dormant}) if the oper connection $\nabla$ is nilpotent of exponent at most $\rk(V)$ (resp. has $p$-curvature zero). Before proceeding further let us recall the following result of \cite[Theorem 3.1.6]{joshi15} which shows that dormant opers of type $q$ exist.

\bthm\label{th:jp316}
Let $Q$ be any vector bundle of rank $q$. Then the triple
$$(F^*(F_*(Q)),\nabla^{can},F^*(F_*(Q))_{\mydot})$$ 
forms a dormant oper of type $q$. Here $\nabla^{can}$ denotes the canonical connection on the vector bundle $F^*(F_*(Q))$.
\ethm

\bigskip

The relationship between opers of type $1$ and the Quot-scheme $\qqquot{Q}$ when $Q$ is a line bundle
was studied in \cite{joshi15}. The main result \cite[Theorem~6.2.1 and Theorem~5.4.1]{joshi15} is the
following 

\bthm \label{opert1}
Assume $p > r(r-1)(r-2)(g-1)$. Let $Q$ be a line bundle of degree $\deg(Q) = -(r-1)(g-1)$. Then the set 
of dormant opers $(V, \nabla, V_\mydot)$ with $V_0/V_1 \isom Q$ is a non-empty, finite set, in bijection
with the set of $k$-rational points $\qqquot{Q}(k)$.

\ethm

\newcommand{\ad}[1]{{\rm ad}(#1)}
\newcommand{\adz}[1]{{\rm ad^0(#1)}}
\newcommand{\omx}{\Omega^1_X}
\newcommand{\ox}{\O_X}
\newcommand{\ext}[3]{Ext^{#1}({#2},{#3})}
\newcommand{\coh}[2]{H^{#1}(X,#2)}
\renewcommand{\hom}{{Hom}}
\newcommand{\br}{\begin{remark}}
\newcommand{\er}{\end{remark}}

We now give an alternative characterization of the underlying bundle $V$ of an oper $(V, \nabla, V_\mydot)$ of type $q$.

\bthm\label{th:oper-con}
Let $X/k$ be a smooth, projective curve over an algebraically closed field $k$. Let $S$ be a vector bundle on $X$. 
Let $\theta$ be a line bundle on $X$. Suppose that the following hypotheses are satisfied:
\benum[label={\bf(\arabic*)}]
\item One has $\theta^2\isom \omx$ i.e. $\theta$ is a theta-characteristic on $X$.
\item $\deg(S)=0$.
\item an isomorphism  $\End(S)=\O_X\oplus \End_0(S)$ of $\ox$-modules.
\item $S$ is stable.
\item $\End_0(S)^* \isom \End_0(S)$.
\eenum
Then for every integer $\ell\geq 2$ there exists a vector bundle $V$ on $X$ satisfying the following 
\benum[label={\bf(\arabic*)}]
\item $\deg(V)=0$
\item $V$ is equipped with a filtration $V=V_0\supset V_1 \supset V_2\supset \cdots \supset V_{\ell}=0$,
\item $\gr i V\isom S \tensor \theta^{2i-(\ell-1)} =: S_i$
\item for all $i=0,1,\ldots, \ell-1$ the extension 
$$0\to \gr{i+1}{V} \to V_i/V_{i+2} \to \gr i V \to 0$$
is the unique non-split extension of this type.
\eenum
Up to an isomorphism, there is only one vector bundle $V$ with these properties.
\ethm

\begin{rem}
If $\rk(Q)$ is not divisible by $p$ then one has an isomorphism $\End(S)=\ox\oplus \End_0(S)$. Further one then also has 
$\End_0(S)^* \isom \End_0(S)$. Hence in particular the hypotheses of the theorem can easily be satisfied for any stable bundle $S$ of 
degree zero and of rank coprime to the characteristic of the ground field. In particular these hypotheses are easily satisfied 
in characteristic zero.
\end{rem}

We begin with the following lemma.

\blem\label{le:con1} 
Let $\theta_i=\theta^{2i-(\ell-1)}$. Then  for all $i\geq 1$,
$$\theta^*_{i-1}\tensor\theta_i\isom \omx.$$
\elem
\bp 
From the definition 
$$\theta^*_{i-1}\tensor\theta_i = \theta^{-(2(i-1)-(\ell-1))+(2i-(\ell-1))}\isom \theta^{2i-2(i-1)}\isom \theta^2 \isom \omx.$$
Hence the claim.
\ep

\blem\label{le:con2}
For $i=0,\ldots,\ell-1$ one has 
$$\ext 1 {S_{i-1}} {S_i} \isom \coh 1 \omx.$$
\elem
\bp 
One has 
\beas
\ext 1 {S_{i-1}} {S_i} &=& \coh 1 {\shom(S_{i-1},S_i)}\\
&=&\coh 1 {S_{i-1}^*\tensor S_i}\\
&=&\coh 1 {S^*\tensor S\tensor \theta_{i-1}^*\tensor\theta_i}\\
&=&\coh 1 {S^*\tensor S\tensor \omx} \text{ [by Lemma~\ref{le:con1}]}.
\eeas
Hence 
\beas 
\ext 1 {S_{i-1}} {S_i} &\isom& \coh 1 {S^*\tensor S\tensor \omx}\\
&=&\coh 1 {(\ox \oplus \End_0(S)) \tensor \omx}\\
&=&\coh 1 \omx \oplus \coh 1 {\End_0(S) \tensor \omx}.
\eeas
By Serre duality 
$$ \coh 1 {\End_0(S) \tensor \omx } \isom \coh 0 {\End_0(S)},$$
and as $S$ is stable of degree zero (by hypothesis) so $\coh 0 {\End_0(S)}=0$. This proves the claim.
\ep

\blem\label{le:con3}
For $i=2,\ldots,\ell-1$ one has
$$\ext 1 {S_{i-2}} {S_i} = 0.$$
\elem
\bp 
One has
\beas 
\ext 1 {S_{i-2}} {S_i} &\isom& \coh 1 {\shom(S_{i-2},S_i)}\\
&=&\coh 1 {S^*\tensor S \tensor \theta_{i-2}^*\tensor \theta_i}\\
&=&\coh 1 {\theta_{i-2}^*\tensor \theta_i \oplus \End_0(S) \tensor \theta_{i-2}^*\tensor\theta_i}.
\eeas
Now 
$$\theta_{i-2}^*\tensor\theta_i\isom \theta^{-(2(i-2)-(\ell-1))+(2i-(\ell-1))}\isom \theta^{-2i+4+2i}\isom {\omx}^{\tensor 2}.$$
Therefore 
$$
\coh 1 {\theta_{i-2}^*\tensor \theta_i} = \coh 1 {{\omx}^{\tensor 2}}=0,
$$
and
$$\coh 1 {\End_0(S) \tensor {\omx}^{\tensor 2}} \isom \coh 0 {\End_0(S) \tensor ({\omx})^{-1}} =0$$
as $S$ is stable, hence there are no global homomorphisms $S \rightarrow S \otimes (\omx)^{-1}$ 
since $\mu(S \otimes (\omx)^{-1}) < \mu(S)$.
This proves the lemma.
\ep

\blem\label{le:con4}
For all $j\geq 3$ one has 
$$H^1(X, \End(S) \tensor\theta^j)=0.$$
\elem
\bp 
This is clear by stability of $S$ and Serre duality.
\ep

\newcommand{\iz}{i_0}           %%% for i_0
\newcommand{\qiz}{S_{\iz}}        % for Q_{i_0}
\newcommand{\qizmo}{S_{\iz-1}}    % for Q_{i_0-1}
\newcommand{\viz}{V_{\iz}}        % for V_{i_0}
\newcommand{\vizpo}{V_{\iz+1}}    % for V_{i_0+1}
\newcommand{\vizmo}{V_{i_0-1}}    % for V_{i_0-1}

\blem\label{le:con5}
For $j_0>i_0+1$ we have
$$Ext^1(\qizmo,S_{j_0-1})=0.$$
\elem

\bp 
We have 
$$Ext^1(\qizmo,S_{j_0-1})\isom H^1(\End(S) \tensor \theta^{2(j_0-\iz)} )$$
and as $j_0 > \iz+1$ so $j_0-\iz>1$. Hence this space is of the form $H^1(\End(S) \tensor \theta^m)$ with $m\geq 3$. 
So we are done by Lemma~\ref{le:con4}.
\ep

\blem\label{le:con6}
Suppose for some $\iz\geq 0$ the bundle $\viz$ has been constructed with asserted properties. Then one has
$$Ext^1(\qizmo,V_j)=0 \qquad \text{ for } j\geq \iz+1.$$
\elem
\bp 
Clearly the claim is true for $j=\ell$ as $V_\ell=0$. We prove the claim by descending induction on $j$. Suppose the lemma is true for some $j_0$ with $\iz+1<j_0\leq \ell$. Then we claim that the assertion is also true for $j_0-1$. As $j_0>\iz+1$, one has the exact sequence
$$0\to V_{j_0}\to V_{j_0-1}\to S_{j_0-1}\to 0.$$
Then applying $Hom(\qizmo,-)$ one gets
$$\to Ext^1(\qizmo,V_{j_0})\to Ext^1(\qizmo,V_{j_0-1}) \to Ext^1(\qizmo,S_{j_0-1})\to 0.$$
By induction hypothesis the first term is zero and by Lemma~\ref{le:con5} the last term is zero. Hence the middle term is zero and the claim is proved.
\ep
Now we are ready to prove Theorem~\ref{th:oper-con}.

\bp[Proof of Theorem~\ref{th:oper-con}]
The construction of the bundle $V$ whose existence is asserted in the theorem is inductive.
We let $V_\ell=0$ and $V_{\ell-1}=S_{\ell-1}$. Now we define $V_{\ell-2}$ as the unique non-split extension 
$$0\to S_{\ell-1}=V_{\ell-1}\to V_{\ell-2} \to S_{\ell-2}\to 0,$$
which is given by the isomorphism of Lemma~\ref{le:con2}:
$$\ext 1 {S_{\ell-2}} {S_{\ell-1}} \isom \coh 1 \omx.$$

Now, if $\ell -3\geq 0$ we define $V_{\ell-3}$ using $V_{\ell-2}$ as follows. Apply $\hom(S_{\ell-3},-)$ to the exact sequence defining $V_{\ell-2}$. This gives
\begin{multline}
0\to \hom(S_{\ell-3},S_{\ell-1})\to \hom(S_{\ell-3},V_{\ell-2})\to \hom(S_{\ell-3},S_{\ell-2})\to \\
\to \ext 1 {S_{\ell-3}} {S_{\ell-1}} \to \ext 1 {S_{\ell-3}} {V_{\ell-2}} \to \ext 1 {S_{\ell-3}} {S_{\ell-2}} \to 0.
\end{multline}
By Lemma~\ref{le:con3} we get $\ext 1 {S_{\ell-3}} {S_{\ell-1}} = 0$, and by Lemma~\ref{le:con2} we get
$$\ext 1 {S_{\ell-3}} {S_{\ell-2}} = \coh 1 \omx.$$
Therefore we get
$$\ext 1 {S_{\ell-3}} {V_{\ell-2}} \isom \coh 1 \omx.$$
So we define $V_{\ell-3}$ as the unique non-split extension given by this isomorphism. Then we have by construction
$$0\to V_{\ell-2} \to V_{\ell-3} \to S_{\ell-3} \to 0.$$
Now we repeat this process to obtain the general construction. Suppose that for 
some $\iz\geq0$, $\viz$ has been constructed with the asserted properties. Then one has an exact sequence
$$0\to \vizpo\to \viz\to \qiz\to 0.$$
Then we claim that there is a unique non-split extension in $Ext^1(\qizmo,\viz)$ which gives $\vizmo$. We proceed as follows: apply $Hom(\qizmo,-)$ to the above short exact sequence to get 
$$\to Ext^1(\qizmo,\vizpo) \to Ext^1(\qizmo,\viz) \to Ext^1(\qizmo,\qiz)\to 0.$$
By Lemma~\ref{le:con2} one as $Ext^1(\qizmo,\qiz)\isom H^1(X,\omx)$, and so to prove our claim it suffices to prove that $Ext^1(\qizmo,\vizpo)=0$. This follows from Lemma~\ref{le:con6}.  Repeating this process   eventually one  gets $V_0$ as the unique non-split extension
 $$0\to V_1\to V_0\to Q_0\to 0.$$
This completes the proof.
\ep

\section{Properties of the Quot-scheme $\qqquot{Q}$ when $q=1$}

The following proposition gives some properties on the Quot-scheme 
$\qqquot{Q}$ when $Q$ is a line bundle. It generalizes the case $r=2$, which was already shown in \cite[Lemma 7.1.3]{joshi15} 
and in the proof of \cite[Theorem 7.1.2]{joshi15}.

\bigskip

We also note that there is a gap in the proof of
\cite[Theorem 7.1.2]{joshi15} and therefore 
we shall give the proof of the following theorem with all the details.

\bigskip

Let $Q$ be a line bundle of 
degree $\deg (Q) = -(r-1)(g-1) + d$ with $d \geq 0$ and let $\sC$ 
be an irreducible component of $\qqquot{Q}$. We recall from Proposition \ref{p:expected-dim}(5) that the expected
dimension of $\sC$ is $rd$.

\begin{defn}
We will say that $\sC$ contains a dormant oper if there exists an effective divisor $D$ of degree $d$ and
a point $[E] \in \qqquot{Q(-D)}$ --- hence $[E]$ corresponds to a dormant oper --- such that 
$[E] \in \sC$ under the natural inclusion
$$ \qqquot{Q(-D)} \subset \qqquot{Q}. $$
\end{defn}

\begin{rem}
We note that the above definition of $\sC$ containing a dormant oper is more restrictive than the
natural one: there exists a $[E] \in \sC$ such that the triple $(F^* E, \nabla^{can}, F^* E_\mydot)$, where
the filtration $F^* E_\mydot$ is the induced filtration from $F^*(F_*(Q))$, is an oper of type $1$.
Note that the latter definition was used in \cite{joshi15} for $r=2$. It can be checked that both
definitions coincide for small $d$.
\end{rem}

\bthm \label{dimopercomponent}
With the above notation and assuming that $p> r(r-1)(r-2)(g-1)$ we have for any irreducible 
component $\sC$ containing a dormant oper
\begin{enumerate}[label={\bf(\arabic*)}]
    \item $\dim (\sC) = rd$.
    \item For a general vector bundle
    $[E] \in \sC$ the map $f_E$ obtained by adjunction
    $$f_E : F^*(E) \rightarrow Q $$
    is surjective.
\end{enumerate}
\ethm

\bp
{\bf (1)} We prove the result by induction on $d$. For $d=0$ this is exactly \cite[Theorem 6.2.1]{joshi15}. 
Assume that the statement holds for an integer $d \geq 0$ and consider a line
bundle $Q$ with $\deg(Q) = -(r-1)(g-1) + (d+1)$. Let $\sC \subset \qqquot{Q}$ be an 
irreducible component containing a dormant oper $[E]$, i.e. by definition
there exists an effective divisor $D$ of degree $d+1$ such that
$E \hookrightarrow F_*(Q(-D)) \subset F_*(Q)$.  From now on the proof goes along the lines
of \cite[Lemma 7.1.3]{joshi15}. For the convenience of the reader we include it here. 

We decompose $D = x + D'$ with $x \in X$
and $D'$ effective of degree $d$. Let $\sC'$ be an irreducible component of 
$\qqquot{Q(-x)} \cap \sC$ containing $[E]$. By induction we have $\dim (\sC') = rd$. Now we 
claim that $\codim_{\sC}{\sC'} \leq r$. To prove this, note that $\sC' \not= \emptyset$. Since 
$\sC$ is an irreducible component of the Quot-scheme, it is equipped with a universal
quotient sheaf $\mathcal{Q}$ over $X \times \sC$.
$$ 0 \longrightarrow \mathcal{E} \longrightarrow p_X^*(F_*(Q)) \longrightarrow 
\mathcal{Q} \longrightarrow 0. $$
We denote by $\mathcal{E}$ the kernel $\ker( p_X^*(F_*(Q)) \longrightarrow 
\mathcal{Q})$. Since $\mathcal{Q}$ and $p_X^*(F_*(Q))$ are $\sC$-flat, $\mathcal{E}$
is also $\sC$-flat and $\forall c \in \sC$ the homomorphism $\mathcal{E}_{|X \times \{ c \}} 
\rightarrow F_*(Q)$ is injective (see e.g. \cite{huybrechts-lehn-book}). Hence, since $F_*(Q)$ is locally free,
$\mathcal{E}_{|X \times \{ c \}}$ is also locally free (since torsion free over a smooth curve)
and by \cite[Lemma 2.1.7]{huybrechts-lehn-book} we conclude that $\mathcal{E}$  is locally free over $X \times \sC$.
Since $p_X^*(F_*(Q)) = (F \times \mathrm{id}_\sC)_* (p_X^* (Q))$ we obtain by adjunction
a  non-zero map $(F \times \mathrm{id}_\sC)^* (\mathcal{E}) \rightarrow p_X^*(Q)$, hence
a non-zero global section $\sigma$ of the rank-$r$ vector bundle 
$$ \mathcal{V} := \shom((F\times {\rm id}_{\sC})^*(\mathcal{E}),p_X^*(Q))$$
over $X \times \sC$. It is clear that $\qqquot{Q(-x)} \cap \sC$ is the zero-scheme
of the section $\sigma_{|\{ x \} \times \sC } \in H^0( \sC, \mathcal{V}_{|\{ x \} \times \sC })$ obtained after restriction.
Hence $\codim_{\sC}{\sC'} \leq r$ and therefore $\dim (\sC) \leq rd + r$. On the other hand, by
the dimension estimates of the Quot-schemes in Proposition \ref{dimensionirrcom}(2)
we have $\dim (\sC) \geq rd+r$. Therefore $\dim (\sC) = rd + r$ and we are done.

\bigskip
{\bf (2)}
First consider the case $d=0$. Then by \cite[Theorem 5.4.1 (1)]{joshi15} any $[E] \in \qqquot{Q}$
is an oper, which implies that $f_E^* : F^*(E) \rightarrow Q$ is surjective.

\bigskip
We now assume that $d >0$. With the previous notation we denote by $\Sigma$ the zero-scheme of 
the global section $\sigma$ of the vector bundle $\mathcal{V}$ and by $p_X$ and $p_{\sC}$ projections onto $X$ and $\sC$ respectively.
Thus one has a diagram
\be 
\xymatrix{
\Sigma \ar[r]& X\times \sC\ar[dl]^{p_X}\ar [dr]_{p_{\sC}} &\\
X  && \sC 
}
\ee

We have the set-theoretical equalities
\begin{eqnarray*}
\Sigma_{|_{\{x\}\times\sC}}&=&\qqquot{Q(-x)}\cap \sC\\
&=&\left\{E\in\sC: f_E: F^*(E)\to Q\ {\rm\ not\ surjective\ at\ }x\right\}.
\end{eqnarray*}

Since $\Sigma$ is closed and $p_{\sC}$ is a proper map, so $\Sigma'=p_{\sC}(\Sigma)$ is a closed subset of the irreducible component  $\sC$. So there are two possibilities:
\indent\benum[label={\bf(\Roman*)}]
\item\label{correction-1} Either  $\Sigma'\neq\sC$, or
\item\label{correction-2} $\Sigma'=\sC$.
\eenum
Now suppose we are in Case~\ref{correction-1}: in this case for a general $[E] \in \sC$, one has $(X\times\{E\})\cap \Sigma=\emptyset$, which is equivalent to surjectivity of $f_E$. So in Case~\ref{correction-1} one has surjectivity of $f_E$ for general $[E]$ as claimed.

So let us suppose we are in Case~\ref{correction-2}. In this case there is at least one irreducible component of $\Sigma$, which we will again denote by $\Sigma$, which surjects onto $\sC$ under $p_{\sC}$. Hence 
$$\dim(\Sigma)\geq\dim(\sC)\geq rd.$$
\newcommand{\pbs}{{\bar p}_X}
Consider the restriction $\pbs:\Sigma\to X$ of $p_X:X\times\sC\to X$ to $\Sigma$. First we assume 
that $\pbs$ is surjective. Since $X$ is smooth and $\Sigma$ is irreducible, by \cite[Chapter~III, Proposition~9.7]{hartshorne-algebraic} $\pbs$ is flat. Hence by \cite[Chapter~III, Corollary~9.6]{hartshorne-algebraic} any irreducible component 
of the fiber $\pbs^{-1}(x)$ has dimension equal to $\dim(\Sigma)-1\geq rd-1$. 
By assumption $\sC$ contains a dormant oper $[E]$ and as $\Sigma'=\sC$, there exists an $x_0\in X$ such that 
$(x_0,[E])\in \Sigma$. Since any irreducible component of $\pbs^{-1}(x_0)$ has dimension $\geq rd-1$, we see that the irreducible component of $\pbs^{-1}(x_0)$ containing $[E]$ has dimension at least $rd-1$ and this irreducible component is contained in $\qqquot{Q(-x_0)}$. But this contradicts part (1) with $d-1 \geq 0$. 

Finally, we need to consider the case when $\pbs$ is not surjective. Since $\Sigma$ is
irreducible, $\pbs(\Sigma) = \{ x_0 \}$ for some point $x_0 \in X$. Hence $\Sigma \subset \pbs^{-1}(x_0)
= \{ x_0 \} \times \sC$, which implies that 
$$ \Sigma \subset \qqquot{Q(-x_0)} \cap \sC. $$
By assumption $\sC$ contains a dormant oper $[E]$, hence there exists an effective divisor $D$ of degree $d \geq 1$
such that $[E] \in \qqquot{Q(-D)}$. Then clearly $x_0 \in D$ and $(x_0, [E]) \in \Sigma$. We then
obtain a contradiction, since $\dim(\Sigma) \geq rd$ and by part (1) the irreducible component
of $\qqquot{Q(-x_0)}$ containing $[E]$ has dimension $r(d-1) < rd$.
\ep

\bigskip

\begin{rem}
The gap in the proof of \cite[Theorem 7.1.2]{joshi15} is related to the fact that  the intersection
$\mathrm{Quot}^{2,0}(F_*(Q(-x))) \cap \sC$ is not necessarily irreducible.
\end{rem}

\section{The Quot-scheme $\qqquot{Q}$ is bigger than expected
when $q \geq 2$}

We assume that $q = \rk(Q) \geq 2$ and that 
$$-(r-q)(g-1) \leq \deg(Q) = -(r-q)(g-1) + d \leq -1,$$
or equivalently,
$$ 0 \leq d \leq (r-q)(g-1) -1.$$
Note that the last inequality implies that $d < (r-1)(g-1) -1$. With this notation we 
have the following result.

\bpro \label{biggerthanexpected}
Let $[E] \in \quot{r,0}{Q}$. Then we have 
$$ \dim \quot{r,0}{Q}  > \mathrm{exp.dim}_{[E]} \ \quot{r,0}{Q} = rd.$$
\epro

\bp
By \cite[Theorem 2.3.1]{joshi15} there exists a line subbundle $L \subset Q$ such that
\beas
\deg L & \geq & \mu(Q) - \left(  1 - \frac{1}{q} \right) (g-1) - \left(1 - \frac{1}{q} \right) \\
       & = & \left( - \frac{r}{q} + 1 \right) (g-1) + \frac{d}{q} - \left(  1 - \frac{1}{q} \right) (g-1) - \left(1 - \frac{1}{q} \right) \\
       & = & - \frac{1}{q}(r-1)(g-1) +  \frac{d}{q} - \left(1 - \frac{1}{q} \right).
\eeas
Note that the expression on the right-hand side is not necessarily an integer. 
We then introduce the quantity $\delta = \deg L + (r-1)(g-1)$. Then clearly 
$$\quot{r,0}{L} \subset \quot{r,0}{Q}$$ 
and by Proposition \ref{dimensionirrcom} (2) $\dim \quot{r,0}{L} \geq r \delta$. Therefore in order to
show the proposition it will be enough to show the inequality $\delta > d$.
We have
\beas
\delta & = & \deg L + (r-1)(g-1), \\
       & \geq & (r-1)(g-1)\left( 1 - \frac{1}{q} \right) + \frac{d}{q}  - \left( 1 - \frac{1}{q} \right) , \\
       & = & \left( 1 - \frac{1}{q} \right) \left[ (r-1)(g-1) -1 \right] + \frac{d}{q}. 
\eeas
Now we observe that the inequality
$$ \left( 1 - \frac{1}{q} \right) \left[ (r-1)(g-1) -1 \right] + \frac{d}{q} > d$$
is equivalent to 
$$ (r-1)(g-1) -1 > d,$$
which holds by assumption on $d$.
\ep

\begin{rem}
In particular if $e(E) = 0$ the Quot-scheme $\quot{r,0}{Q}$ has some positive-dimensional
components, as already observed in the previous section, see also Theorem \ref{p:quot-product-lemma}.
\end{rem}

\bigskip

\section{A characterization of dormant opers of type $q$}

\blem\label{l:oper-flag-lemma}
Let $Q$ be any vector bundle of rank $q$ and $\deg(Q)= -q(\ell -1)(g-1)$. We denote by $\trip V$ the 
oper $(F^*(F_*(Q)),\nabla^{can},F^*(F_*(Q))_{\mydot})$ introduced in Theorem \ref{th:jp316}. 
Let $[W'] \in\qquot$ and let $W=F^*(W') \subset V$. Then the triple $(W, \nabla^{can}, W_\mydot)$, where 
$\nabla^{can}$ is the canonical connection on $F^*(W')$ and $W_\mydot$ the induced flag defined by
$W_i = W \cap V_i$, is an oper of type $q$ if and only if 
$$W_\ell=W\cap V_\ell= \{ 0\}. $$
\elem

%\medskip

\bp 
Suppose that the triple  $(W, \nabla^{can}, W_\mydot)$ is an oper of type $q$. By \cite[Lemma 3.4.2(i)]{joshi15} 
we have an inclusion $W_0/W_1 \subset Q$ and since both vector bundles have same degree and rank, they are isomorphic.
It also follows from \cite[Lemma 3.4.2(ii)]{joshi15} that $\rk(W_\ell) = 0$, hence $W_\ell = \{ 0 \}$.
\bigskip

Conversely, assume that $W_\ell= \{ 0 \}$. Then consider the sequence $r_i=\rk(W_i/W_{i+1})$. So
by \cite[Lemma~3.4.2(ii)]{joshi15} we have the inequalities 
$$q\geq r_0\geq r_1\geq\cdots\geq r_\ell=0.$$
As $\sum_{i=0}^{\ell-1}r_i=q\ell$ it follows that all $r_i=q$. Thus $\rk(W_i/W_{i+1})=r_i=q$ for all $i$.

Now consider the increasing sequence of degrees $\deg(V_i/V_{i+1})=\deg(Q)+iq(2g-2)$ for $i=0,\ldots,\ell-1$.
Since we have an injective map $W_i/W_{i+1}\into V_i/V_{i+1}$ and both the bundles are of the same rank, we obtain 
$\deg(W_i/W_{i+1})\leq \deg(V_i/V_{i+1})$ and one has
$$0=\sum_{i=0}^{\ell-1}\deg(W_i/W_{i+1}) \leq \sum_{i=0}^{\ell-1}\deg(V_i/V_{i+1}).$$
This gives
$$0=\sum_{i=0}^{\ell-1}\deg(W_i/W_{i+1}) \leq \sum_{i=0}^{\ell-1}(\deg(Q)+iq(2g-2)).$$
The last sum evaluates to
$$0=\sum_{i=0}^{\ell-1}\deg(W_i/W_{i+1}) \leq q(g-1)(\ell(\ell-1)-\ell(\ell-1))=0.$$
Thus one sees that $\deg(W_i/W_{i+1})=\deg(V_i/V_{i+1})$. Therefore the above injective maps are isomorphisms:
$W_i/W_{i+1} \isom  V_i/V_{i+1}$ for $i=0,\ldots,\ell-1$. This proves the assertion.
\ep

One important consequence of this lemma is the following fundamental result. 

\bthm\label{th:dormant-q-opers}
Let $Q$ be a semi-stable vector bundle with $\deg(Q)= - q(\ell-1)(g-1)$ and $\rk(Q) = q$ for some 
integer $\ell \geq 2$. 
We put $r = \ell q$ and use the notation of Lemma \ref{l:oper-flag-lemma}.
\benum[label={\bf(\arabic*)}]
\item For every $[W']\in\qquot$ satisfying 
$$ W_\ell = W \cap V_\ell = \{0\}$$ 
the triple $(W = F^*(W'), \nabla^{can}, W_\mydot)$
is a dormant oper of type $q$.
\item If $p>r(r-1)(g-1)$ then conversely every dormant oper $\trip{W}$ of type $q$ with $W_0 / W_1 \isom Q$  is of the 
form $(W = F^*(W'), \nabla^{can}, W_\mydot)$ for some $W'\in\qquot$.
\eenum
\ethm 

\bp
The first assertion is immediate from Lemma \ref{l:oper-flag-lemma}. 

\smallskip
Now suppose $\trip W$ is a dormant oper of type $q$. Then by Cartier's Theorem (see \cite{katz70}) 
$W=F^*(W')$ for some vector bundle $W'$ and as $W=F^*(W')\twoheadrightarrow W_0/W_1=Q$ one gets by adjunction 
a non-zero map $W' \to  F_*(Q)$. Let us first show that $W'$ is semi-stable (of degree zero). 
Suppose this is not the case. Then there exists a subbundle $W''\into W'$ 
for which $\mu(W'')>0$. Then in fact one has $\mu(W'')\geq \frac{1}{r-1}$. But then 
$F^*(W'')\into W$ such that $\mu(F^*(W''))=p\mu(W'')\geq \frac{p}{r-1}$. On the other hand as $W$ 
carries the structure of an oper of type $q$, and one has $\mu(W_{\ell-1})=q(2(\ell-1)-(\ell-1))(g-1)=q(\ell-1)(g-1)$ 
and this is the destabilizing subbundle of largest degree, so $\frac{p}{r-1}\leq \mu(F^*(W''))\leq q(\ell-1)(g-1)$. 
Hence
$$p<(r-1)q(\ell-1)(g-1)=(r-1)(r-q)(g-1).$$ On the other hand we have assumed that $p>r(r-1)(g-1)$. Therefore we have arrived at a contradiction. Thus the vector bundle $W'$ is semi-stable. 

\bigskip

Now if $W'$ does not map injectively into $\fs{Q}$, then the image is a subsheaf of some degree $d\geq 1$ and 
rank $\leq r-1$ and hence it has slope $\geq\frac{1}{r-1}$. Again by \cite[Proposition 4.2.1]{joshi15} one sees that $\fs{Q}$ does not have any subsheaves of suitably positive degree and of rank $\leq r-1$. So $W' \into \fs{Q}$.
\ep

\bigskip

\bpro
Using the above notation, the set of points $[W'] \in \qquot$ satisfying $W_\ell = \{ 0 \}$ is an open subset.
\epro

\bp
The Quot-scheme $\qquot$ comes equipped with a universal quotient $\sG$ defined over the product
$X\times \qquot$ which is flat over $\qquot$:
\be 
0\to \sW' \to p_X^* (F_*(Q)) \to\sG \to 0.
\ee 
Since $p_X^* (F_*(Q))$ is locally free, it is flat over $\qquot$ and therefore for any point $y \in \qquot$ the
homomorphism $\sW'|X \times \{ y \} \rightarrow F_*(Q)$ is injective. Thus $\sW'|X \times \{ y \}$ is locally free for any 
$y$, which in turn implies that $\sW'$ is locally free. Hence, with the notation of Lemma \ref{l:oper-flag-lemma} we obtain a homomorphism between locally free sheaves over 
$X \times \qquot$
$$ \Phi:  F^*(\sW') \rightarrow p_X^* (V/V_\ell). $$
By Lemma \ref{l:oper-flag-lemma} the condition $W_\ell = \{0 \}$, with $W = F^* (\sW'|X \times \{ y \})$,
is equivalent to  $\Phi_{X \times \{ y \} }$ being an isomorphism. But the set of points $y \in \qquot$ 
where $\Phi_{X \times \{ y \}}$ is not an isomorphism is clearly a closed subset. 
\ep

\bigskip

We shall denote this open subscheme $$\Omega(Q) \subset \qquot. $$ 
We recall that by Theorem \ref{th:dormant-q-opers} the open 
subscheme $\Omega(Q)$ 
parameterizes dormant opers $\trip W$ of type $q$ with $W_0 / W_1 \isom Q$.

\bigskip

\begin{rem}
We observe that the Quot-scheme $\qquot$ has expected 
dimension $0$ if $\rk(Q) = q$, $\deg(Q)= - q(\ell-1)(g-1)$ and $r = q \ell$, but in the case $q \geq 2$ we have 
shown in Proposition \ref{biggerthanexpected} that $\qquot$ always contains a Quot-scheme 
$\mathrm{Quot}^{r,0}(F_*(L))$ for some line subbundle $L \subset Q$, which has dimension $> 0$.
Clearly, the Quot-scheme $\mathrm{Quot}^{r,0}(F_*(L))$ is not contained in the open subscheme 
$\Omega(Q)$ and does not correspond to dormant opers.
\end{rem}

\bcon
Let $Q$ be a semi-stable vector bundle with $\deg(Q)= - q(\ell-1)(g-1)$ and $\rk(Q) = q$ for some
integer $\ell \geq 2$. Put $r = q\ell$. Then, for $p > r(r-1)(r-2)(g-1)$ 
$$ \Omega(Q) \text{ is non-empty and of dimension} \ 0. $$
\econ

\smallskip

We now list some evidence for this conjecture:

\begin{enumerate}
    \item If $q= 1$, the conjecture is true. This is shown in \cite{joshi15} --- see Theorem \ref{opert1}.
    Note that in this case we have equality $\Omega(Q) = \qquot$.
    \item The conjecture is true for a decomposable bundle of the form $Q = \oplus_{i=1}^q L_i$, where the 
    $L_i$ are $q$ distinct line bundles of degree $\deg(L_i) = - (\ell -1)(g-1)$. This is shown in 
    Theorem \ref{p:quot-product-lemma}. We note that in this case points of $\Omega(Q)$ correspond
    to direct sums of $q$ opers of type $1$.
    
\end{enumerate}

\newcommand{\qmod}[1]{{\mathscr M}_X(q,{#1})}

\section{Dimension of Frobenius instability loci}

\newcommand{\sM}{{\mathcal M}}
\newcommand{\sJ}{\mathcal{J}}
\newcommand{\sU}{\mathcal{U}}
\newcommand{\umod}{{\sM}(r)}

Throughout this section we assume that $p > r(r-1)(r-2)(g-1)$.

\bigskip

Let $\umod$ be the coarse moduli space of semi-stable vector bundles of rank $r$ and degree $0$ over $X$.
Let $\sJ(r)\subset \umod$ be the closed subscheme of $\umod$ parameterizing semi-stable vector bundles $E$ such that
$F^*(E)$ is not semi-stable. The closed subscheme $\sJ(r)$ will be referred to as the Frobenius instability locus. 
The purpose of this section is to make a conjecture on the dimension of $\sJ(r)$. Before stating the conjecture
we recall some facts on the structure of $\sJ(r)$.

\bigskip

Let $\sM(q,-1)$ be the moduli space of semi-stable bundles of rank $q$ and degree $-1$ over $X$. 
Note that as we are in the coprime case any semi-stable bundle of rank $q$ and degree $-1$ is also stable.
Let $\sU$ be the universal bundle over $X \times \sM(q,-1)$. Consider the relative Quot-scheme
$$\pi:  \mathrm{Quot}^{r,0}((F \times \mathrm{id}_{\sM(q,-1)})_* \sU) \longrightarrow \sM(q,-1) $$
with fiber $\pi^{-1}(Q) = \qqquot{Q}$ over a stable bundle $Q \in \sM(q,-1)$.

\bigskip

For $q = 1, \ldots, r-1$ we denote by $\sJ_q \subset \sJ(r)$ the closure of the forgetful map
$$ \alpha_q : \mathrm{Quot}^{r,0}((F \times \mathrm{id}_{\sM(q,-1)})_* \sU)  \longrightarrow \umod, \ \ 
[E \subset F_*(Q)] \rightarrow E.$$
Then by \cite[Theorem 4.4.1]{joshi15} we have the inclusions
$$ \sJ^s(r) \subset \bigcup_{q=1}^{r-1} \sJ_q \subset \sJ(r), $$
where $\sJ^s(r)$ denotes the subset of stable vector bundles $E$ such that $F^*(E)$ is not 
semi-stable.

\bigskip

In order to compute the dimension of $\sJ_q$ we need to know the following
dimensions
\begin{enumerate}
    \item $\dim \qqquot{Q}$ for general $Q \in \sM(q,-1)$.
    \item $\dim \alpha_q^{-1}(E)$ for general $E \in \sJ_q$.
\end{enumerate}

Unfortunately we only know $\dim \sC$ for certain irreducible components $\sC \subset
 \qqquot{Q}$ when $q=1$.
 
\bigskip

We therefore focus on the case $q=1$. Let $Q$ be a line bundle of degree $-1$ and let 
$\sC$ be an irreducible component of $\qqquot{Q}$ containing a dormant oper. 
Then by Theorem \ref{dimopercomponent} $\dim \sC = r((r-1)(g-1) -1)$ and for a general
vector bundle $[E] \in \sC$ the map $f_E : F^*(E) \rightarrow Q$ obtained by 
adjunction is surjective. Let us denote the kernel of $f_E$ by $S = \ker f_E$. Note that
$\deg(S) = 1$ for general $[E]$.

\bigskip

With this notation we make the following 

\bcon \label{condimil}
For a general vector bundle $E \in \sC$ the bundle $S$ is semi-stable.
\econ 

As a consequence of this conjecture we obtain that for a general $E \in \sC$ the
filtration $0 \subset S \subset F^*(E)$ is the Harder-Narasimhan filtration of 
$F^*(E)$. Therefore the quotient $F^*(E) \twoheadrightarrow Q = F^*(E)/S$ is
uniquely determined and $\dim \Hom(E, F_*(Q)) = 1$. This means that the forgetful
map $\alpha_1$ is generically injective on components $\sC$ containing dormant opers.

\bigskip

Hence if Conjecture \ref{condimil} holds, then we have the following formula

$$ \dim \sJ(r) \geq \dim \sJ_1 \geq \dim \sM(q,-1) + \dim \sC = (r^2 - r + 1)(g-1) - (r-1). $$

\begin{rem}
Conjecture \ref{condimil} holds trivially for $r=2$ and the lower bound of the
dimension of $\sJ(2)$ was already worked out in \cite[Theorem 7.1.2]{joshi15}.
\end{rem}

\bigskip

Finally, we conjecture the following structure of $\sJ_q$ for any $q=1, \ldots , r-1$.

\bcon
The Harder-Narasimhan filtration of $F^*(E)$ of a general bundle $E \in \sJ_q$ is
minimal, i.e., of the form
$$ 0 \subset S \subset F^*(E), $$
where $\deg(S) = 1$ and $\rk(S) = r-q$.
\econ

\bigskip

\begin{rem}
The above two conjectures hold in the case $p=2,r=2$ by \cite{jrxy06} and for $g=2, r=3, p=3$. This is worked out in \cite{li16}.
\end{rem}

%\bibstyle{alpha}
%\bibliographystyle{alpha}

%\bibstyle{epigaref}
%\bibliographystyle{epigaref}

\end{document}